\documentclass{amsart}

\usepackage{graphicx,tikz,amsmath}

\newtheorem{theorem}{Theorem}[section]
\newtheorem{lemma}[theorem]{Lemma}
\newtheorem{prop}[theorem]{Proposition}

\DeclareMathOperator{\olw}{olw}
\DeclareMathOperator{\sz}{Sz}

\author{Csaba Bir\'o}
\email{csaba.biro@louisville.edu}
\address{Department of Mathematics, University of Louisville, Louisville, KY 40292}

\author{Israel R. Curbelo}
\email{icurbelo@kean.edu}
\address{Department of Mathematics, University of Louisville, Louisville, KY 40292}
\address{Current address: Department of Mathematical Sciences, Kean University, Union, NJ 07083}

\title[On-line chain partitioning of $d$-dimensional posets]{Improved lower bounds on the on-line chain partitioning of posets of bounded dimension}

\begin{document}

\begin{abstract}
        An on-line chain partitioning algorithm receives a poset, one element at a time, and irrevocably assigns the element to one of the chains. Over 30 years ago, Szemer\'edi proved that any on-line algorithm could be forced to use $\binom{w+1}{2}$ chains to partition a poset of width $w$. The maximum number of chains that can be forced on any on-line algorithm remains unknown. In a survey paper by Bosek et al., it is shown that Szemer\'edi's argument could be improved to obtain a lower bound almost twice as good. Variants of the problem were considered where the class is restricted to posets of bounded dimension or where the poset is presented via a realizer of size $d$. In this paper, we prove two results. First, we prove that any on-line algorithm can be forced to use $(2-o(1))\binom{w+1}{2}$ chains to partition a $2$-dimensional poset of width $w$. Second, we prove that any on-line algorithm can be forced to use $(2-\frac{1}{d-1}-o(1))\binom{w+1}{2}$ chains to partition a poset of width $w$ presented via a realizer of size $d$.

\end{abstract}

\maketitle

\section{Introduction}



We consider each problem as a two-player coloring game between Beth and Anna. In this game, Beth is the builder who constructs a poset one point at a time and Anna is the algorithm who constructs a chain partition. During round $i$, Beth introduces a new point $x_i$ to the poset and describes the subposet induced by the elements $\{x_1,\ldots,x_i\}$. Anna responds by assigning $x_i$ to one of the chains in the chain partition. To avoid confusion, we refer to the chains in the partition as colors. 

The on-line width $\text{olw}(w)$ of the class of posets of width at most $w$ is the largest integer $k$ for which there exists a strategy for Beth that forces Anna to use $k$ colors on a poset of width at most $w$. 
Clearly, $\olw(1)=1$. Kierstead \cite{kie-81} showed that $5\leq \olw(2)\leq6$, and Felsner \cite{fel-97} later showed that $\olw(w)\leq 5$, solving the problem for $w=2$. The exact value of $\olw(w)$ remains unknown for $w\geq 3$. Kierstead \cite{kie-81} was the first to prove that $\olw(w)$ was bounded. The upper bound has since been improved several times with the most recent coming in the year $2021$ from Bosek and Krawcyk \cite{bos-kra-21} where they prove that $\olw(w)\leq w^{O(\log \log w)}$. On the other hand, a strategy by Szemer\'edi \cite{kie-86} proved that $\olw(w)\geq \binom{w+1}{2}$. Szemer\'edi's strategy was later improved in \cite{bos-fel-12} to show that $\text{olw}(w)\geq (2-o(1))\binom{w+1}{2}$. 

In this paper, we show that any poset resulting from the improved strategy is $2$-dimensional. Furthermore, if Beth is required to present the poset along with a realizer of size $d$, then we can achieve a lower bound which is arbitrarily close to that of the general problem for sufficiently large $d$. 



\subsection{Background and Results}

A set $R=\{L_1,\ldots,L_t\}$ of linear extensions of a poset $(X,P)$ is called a realizer of $(X,P)$ if $x<y$ in $P$ if and only if $x<y$ in $L_i$ for $i\in\{1,\ldots,t\}$. The dimension of $(X,P)$ is then defined as the least integer $d$ for which $(X,P)$ has a realizer of cardinality $d$. In this paper, we focus on variants of the problem where not only is the width of the poset restricted but also the dimension of the poset. We refer the reader to the survey paper \cite{bos-fel-12} for an overview of the different variants that have branched from the main problem. Let $\olw(w,d)$ be the largest integer $k$ for which there exists a strategy for Beth that forces Anna to use $k$ colors on a poset of width at most $w$ and dimension at most $d$. The analysis of the on-line chain partition game restricted to $d$-dimensional posets appears to be as hard as the general problem and no better upper bound is known for this class (even for $d=2$). In \cite{bos-fel-12}, a proof that $\text{olw}(w,2)\geq \binom{w+1}{2}$ is provided. For our first contribution, we prove the following theorem.

\begin{theorem}
Let $w$ and $d$ be integers such that $w\leq 1$ and $d\geq 2$. Then $\olw(w,d)\geq (2-o(1))\binom{w+1}{2}$.
\end{theorem}

For the main result of this paper, we consider the variant of the problem first analyzed by Kierstead, McNulty, and Trotter \cite{kie-mcn-84} in which Beth introduces a $d$-dimensional poset via its embedding in $\mathbb{R}^d$ or equivalently, by providing on-line a realizer of cardinality $d$. Let $\olw_R(w,d)$ be the largest integer $k$ for which there exists a strategy for Beth that forces Anna to use $k$ chains on a poset of width $w$ introduced on-line via a realizer of cardinality $d$. Kierstead, McNulty and Trotter \cite{kie-mcn-84} proved that $\olw_R(w,d)\leq\binom{w+1}{2}^{d-1}$. For our main contribution, we prove the following theorem.

\begin{theorem}
Let $w$ and $d$ be positive integers. Then $\olw_R(w,d)\geq (2-\frac{1}{d-1}-o(1))\binom{w+1}{2}$.
\end{theorem}

In Section $2$, we introduce substrategies and generalize previous results. In Sections $3$ and $4$ we prove Theorem 1.1 and 1.2 respectively. First, we introduce some notation and terminology.

\subsection{Notation}
Let $(X,P)$ be a poset. Let $U$ and $V$ be disjoint subsets of $X$. We say that $U < V$ if for any point $u\in U$ and any point $v\in V$, $u<v$. We say that $U$ and $V$ are completely comparable if for any point $u\in U$ and any point $v\in V$, $u$ and $v$ are comparable. Similarly, we say that $U$ and $V$ are completely incomparable if for any point $u\in U$ and any point $v\in V$, $u$ and $v$ are incomparable. Let $T(w)$ be a strategy. Suppose $T(w)$ is played resulting in the poset $(X,P)$. If $Y\subseteq X$, we use $\sigma(Y)$ to denote the set of colors used on $Y$ and $\|Y\|=|\sigma(Y)|$ The dual $P*$ of a partial order $P$ is a partial order on the same set at $P$ such that $x<y$ in $P^*$ if and only if $y<x$ in $P$. Hence, we define the dual strategy $T^*(w)$ as a strategy such that whenever $Beth$ introduces a new point $x$, $x<y$ in $P$ for any previously introduced element $y$ if and only if it would have been the case that $y<x$ in $P^*$ playing $A^*(w)$. Lastly, each strategy in this paper is defined in multiple stages. We reserve $S_i$ to denote the set of points introduced in Stage $i$.
\section{Preliminaries}

\begin{figure}
\centering
\begin{tikzpicture}[scale=0.5]

\def\bot{-2}
\def\top{7}
\def\lef{-1}
\def\rig{21}

\draw[fill=black] (0,0) circle (2pt) node[above left] {$\scriptstyle x_1^1$};
\draw[fill=black] (1,1) circle (2pt) node[above left] {$\scriptstyle x_2^2$};
\draw[fill=black] (3,3) circle (2pt) node[above left] {$\scriptstyle x_4^3$};
\draw[fill=black] (6,6) circle (2pt) node[above left] {$\scriptstyle x_7^4$};

\draw[fill=black] (20,2) circle (2pt) node[above right] {$\scriptstyle x_3^2$};
\draw[fill=black] (18,4) circle (2pt) node[above right] {$\scriptstyle x_5^1$};
\draw[fill=black] (17,5) circle (2pt) node[above right] {$\scriptstyle x_6^3$};

\draw[dashed] (\lef,-1) -- (11.5,-1) node[below]{$\sz(3)$} -- (\rig,-1);
\draw[dashed] (7,\bot) -- (7,\top);
\draw[dashed] (16,\bot) -- (16,\top);

\draw (7,-1) node[above left]{$C_4$};
\draw (16,-1) node[above right]{$G_4$};

\draw[thick] (0,0) -- (6,6);
\draw[thick] (1,1) -- (20,2);
\draw[thick] (3,3) -- (18,4);
\draw[thick] (3,3) -- (17,5);

\draw (\lef,\bot) -> (\rig,\bot) node[right] {$L_\alpha$};
\draw (\lef,\bot) -> (\lef,\top) node[above] {$L_\beta$};
\end{tikzpicture}
\caption{$\sz(4)$ embedded in $\mathbb{R}^2$. Subscripts denote the on-line order and superscripts denote the color assigned. }
\label{fig:n1}
\end{figure}
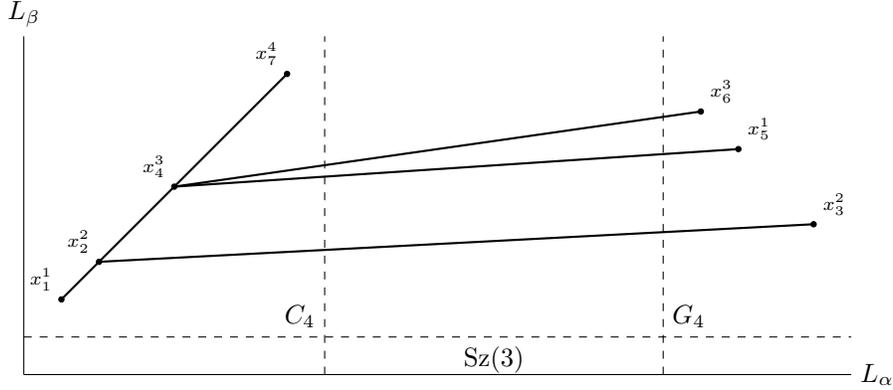

Let us recall Szemer\'edi's original strategy $\sz(w)$ which consists of two stages. In Stage 1, Beth builds a chain $C_w$ of size $w$ along with an antichain $G_w$ of size at most $w-1$. In Stage 2, Beth recursively plays $\sz(w-1)$ so that each new element is less than every element of $G_w$ and incomparable to every element of $C_w$. Szemer\'edi's strategy always results in a poset $(X,P)$ consisting of a chain $C_i$ of size $i$ such that every element of $C_i$ has a distinct color and an antichain $G_i$ of size at most $i-1$ for each $i\in\{1,\ldots, w\}$. We refer to the chains $C_1,\ldots,C_w$ as \emph{rainbow chains} where the index indicates the size of the chain. Moreover, elements from distinct rainbow chains are incomparable and the set of minimal elements of $(X,P)$ are exactly the minimal elements of $C_1,\ldots,C_w$. Figure \ref{fig:n1} shows a poset constructed by $\sz(4)$ embedded in $\mathbb{R}^2$. The $x$ and $y$ coordinates provide linear orders $L_\alpha$ and $L_\beta$ respectively such that $L_\alpha\cap L_\beta=P$. Notice that this embedding strategy will always conclude with $C_w$ at the bottom of $L_\alpha$. 

\subsection{Substrategies}
In order to prove our results, we need to be able to guarantee that the strategy concludes with $C_k$ at the bottom of $L_\alpha$ for any $k\in\{1,\ldots,w\}$. To guarantee this, we must present two strategies $L_\alpha(k,w)$ and $L_\beta(k,w)$ which construct the linear orders $L_\alpha$ and $L_\beta$ respectively. The two strategies are played simultaneously to construct the partial order $L_\alpha\cap L_\beta$, however, in order to acheive our desired condition, we do in fact need to treat them as independent strategies. We prove the existence of such strategies in Lemma \ref{lem2}, but first we prove the following claim.

\begin{lemma}\label{lem1}
    Let $w$ be a positive integer. Then there exist strategies $L_\alpha(w)$ and $L_\beta(w)$ that constructs a realizer $\{L_\alpha,L_\beta \}$ so that $L_\alpha\cap L_\beta=C_w\cup G_w$ where $C_w$ is a rainbow chain of size $w$, $G_w$ is an antichain of size at most $w-1$, $C_w\cap G_w = \emptyset$ and $C_w<G_w$ in $L_\beta$.
\end{lemma}

\begin{proof}
    We argue by induction on the positive integer $w$. If $w=1$, Beth simply presents a single point. Suppose $w>1$. By the induction hypothesis, there exist strategies that constructs a realizer $\{L_\alpha,L_\beta \}$ so that $L_\alpha\cap L_\beta=C_{w-1}\cup G_{w-1}$ where $C_{w-1}$ is a rainbow chain of size $w-1$, $G_{w-1}$ is an antichain of size at most $w-2$, and $C_{w-1}<G_{w-1}$ in $L_\beta$. Beth plays the strategies and then introduces a new point $x$ at the top of $L_\alpha$ and in between $C_{w-1}$ and $G_{w-1}$ in $L_\beta$. If $x$ is assigned a new color, then $C_w=C_{w-1}\cup\{x\}$ and $G_w=G_{w-1}$. Otherwise, we add $x$ to $G_{w-1}$ and repeat. If $|G_{w-1}|=w-1$, then Anna must use a new color in the next round. Thus the strategy terminates with $|C_w|=w$ and $|G_w|\leq w-1$.
\end{proof}

Note that the strategies from Lemma \ref{lem1} can be defined independently in a constructive manner. Each round, $L_\alpha(w)$ simply places a new element at the top of $L_\alpha$, and $L_\beta(w)$ traverses up $L_\beta$ and inserts a new element immediately under the first repeated color. It is easy to verify that played together, they indeed satisfy the conditions of Lemma \ref{lem1}. We now use Lemma \ref{lem1} to prove the following. 

\begin{lemma}\label{lem2}
Let $w$ and $k$ be positive integers such that $k\leq w$. Then for any on-line algorithm there exist strategies $L_\alpha(k,w)$ and $L_\beta(k,w)$ for Beth which construct the same poset as $\sz(w)$ via a realizer $\{L_\alpha,L_\beta\}$ so that if $C_k$ is the rainbow chain of size $k$, then if $u\in C_k$ and $v\notin C_k$, $u<v$ in $L_\alpha$. 
\end{lemma}

\begin{proof}
We argue by induction on the positive integers $w$ and $k$ with $k\leq w$. If $w=k=1$, Beth simply presents a single point. Suppose $w>1$ and $k< w$. The strategies consists of two stages.

By Lemma \ref{lem1}, there exist strategies $L'_\alpha(w)$ and $L'_\beta(w)$ which construct a realizer $\{L'_\alpha,L'_\beta \}$ so that $L'_\alpha\cap L'_\beta=C_w\cup G_w$ where $C_w$ is a rainbow chain of size $w$, $G_w$ is an antichain of size at most $w-1$, and $C_w<G_w$ in $L'_\beta$. By the induction hypothesis, there exist strategies $L''_\alpha(k,w-1)$ and $L''_\beta(k,w-1)$ for Beth which construct the same poset as $\sz(w-1)$ via a realizer $\{L''_\alpha,L''_\beta\}$ so that the rainbow chain $C_k$ of size $k$ is at the bottom of $L''_\alpha$. Then Beth plays $L'_\alpha(w)$ and $L'_\beta(w)$ in Stage 1 followed by $L''_\alpha(k,w-1)$ and $L''_\beta(k,w-1)$ in Stage 2 in $L_\alpha$ and $L_\beta$ respectively so that $S_2<S_1$ in $L_\alpha$ and $C_w<S_2<G_w$ in $L_\beta$.

If $k=w$, then the argument is nearly identical except that Beth uses the induction hypothesis strategies for $k=w-1$, and the roles of $L_\alpha$ and $L_\beta$ are switched.
\end{proof}

\begin{figure}
\centering
\resizebox{4.75in}{!}{\tikzstyle{top_left}=[rectangle,draw,fill=white,minimum size=1.4em,font={$L_\alpha (w-1,w-1)$}]
\tikzstyle{top_right}=[rectangle,draw,fill=white,minimum size=1.4em,font={$L_\alpha (k,w-1)$}]
\tikzstyle{bottom_left}=[rectangle,draw,fill=white,minimum size=1.4em,font={$L_\beta (w-1,w-1)$}]
\tikzstyle{bottom_right}=[rectangle,draw,fill=white,minimum size=1.4em,font={$L_\beta (k,w-1)$}]
\begin{tikzpicture}

\draw (-8,0) -- (8,0);
\draw (0,-2) -- (0,2);
\draw (-8,-2) rectangle (8,2);
\node[above] at (-4,2) {$k=w$};
\node[above] at (4,2) {$k<w$};
\node[left] at (-8,1) {$L_\alpha$};
\node[left] at (-8,-1) {$L_\beta$};

\draw[fill=black] (-0.5,-1) circle (1pt) node[below] {$x_N^w$};
\draw[fill=black] (-1.25,-1) circle (1pt) node[below] {$x_{N-1}^{w-1}$};
\draw[fill=black] (-2.5,-1) circle (1pt) node[below] {$x_4^2$};
\draw[fill=black] (-3,-1) circle (1pt) node[below] {$x_3^2$};
\draw[fill=black] (-3.5,-1) circle (1pt) node[below] {$x_2^1$};
\draw[fill=black] (-4,-1) circle (1pt) node[below] {$x_1^1$};
\draw (-0.5,-1) -- (-1.5,-1);
\draw[dotted] (-1.5,-1) -- (-2.25,-1);
\draw (-2.25,-1) -- (-4.5,-1) node[draw,left] {$L_\beta (w-1,w-1)$};

\draw[fill=black] (7.5,1) circle (1pt) node[below] {$x_N^w$};
\draw[fill=black] (6.75,1) circle (1pt) node[below] {$x_{N-1}^{w-1}$};
\draw[fill=black] (5,1) circle (1pt) node[below] {$x_4^2$};
\draw[fill=black] (4.5,1) circle (1pt) node[below] {$x_3^2$};
\draw[fill=black] (4,1) circle (1pt) node[below] {$x_2^1$};
\draw[fill=black] (3.5,1) circle (1pt) node[below] {$x_1^1$};
\draw (7.5,1) -- (6.5,1);
\draw[dotted] (6.5,1) -- (5.25,1);
\draw (5.25,1) -- (3,1) node[draw,left] {$L_\alpha (k,w-1)$};

\draw (-7.5,1) -- (-6.75,1);
\draw[dotted] (-6.75,1) -- (-6.25,1);
\draw (-6.25,1) -- (-1.75,1);
\draw[dotted] (-1.75,1) -- (-1.25,1);
\draw (-1.25,1) -- (-0.5,1);
\draw[fill=black] (-7.5,1) circle (1pt) node[below] {$x_1^1$};
\draw[fill=black] (-7,1) circle (1pt) node[below] {$x_3^2$};
\draw[fill=black] (-6,1) circle (1pt) node[below] {$x_N^w$};
\node at (-4,1) [top_left] {};
\draw[fill=black] (-0.5,1) circle (1pt) node[below] {$x_2^1$};
\draw[fill=black] (-1,1) circle (1pt) node[below] {$x_4^2$};
\draw[fill=black] (-2,1) circle (1pt) node[below] {$x_{N-1}^{w-1}$};

\draw (7.5,-1) -- (6.75,-1); 
\draw[dotted] (6.75,-1) -- (5.85,-1);
\draw (5.85,-1) -- (2.25,-1);
\draw[dotted] (2.25, -1) -- (1.25,-1);
\draw (1.25,-1) -- (0.5,-1);
\draw[fill=black] (0.5,-1) circle (1pt) node[below] {$x_1^1$};
\draw[fill=black] (1,-1) circle (1pt) node[below] {$x_3^2$};
\draw[fill=black] (2.5,-1) circle (1pt) node[below] {$x_N^w$};
\node at (4,-1) [bottom_right] {};
\draw[fill=black] (7.5,-1) circle (1pt) node[below] {$x_2^1$};
\draw[fill=black] (7,-1) circle (1pt) node[below] {$x_4^2$};
\draw[fill=black] (5.6,-1) circle (1pt) node[below] {$x_{N-1}^{w-1}$};

\end{tikzpicture}}
\caption{The resulting $L_\alpha$ and $L_\beta$ on a greedy algorithm.}
\label{fig:n0}
\end{figure}
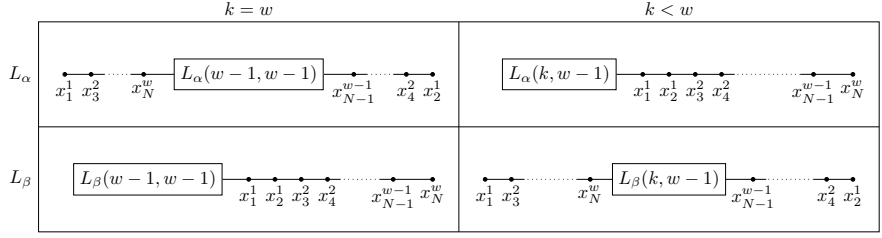

The construction of $L_\alpha$ and $L_\beta$ can be seen in Figure \ref{fig:n0}. The strategies are defined recursively where $L_\alpha(k,w)$ continues to play $L_\alpha(i-1)$ under $L_\alpha(i)$ until $i=k+1$ in which case $L_\beta(k)$ is played under $L_\alpha(k+1)$. This switch when $k=w$ guarantees that the rainbow chain $C_k$ is at the bottom of $L_\alpha$. 

\subsection{Generalizing the improved strategy}

The improved strategy $\sz'(w)$ presented in \cite{bos-fel-12} involved three stages where Beth played $\sz(w)$ in Stage 1 followed by $\sz^*(w)$ completely under $S_1$ in Stage 2. Let $C_1,\ldots,C_w$ be the rainbow chains from $\sz(w)$ and $D_1,\ldots,D_w$ be the rainbow chains from $\sz^*(w)$. Then, the proof relied on the fact that there exists an integer $t$ such that $\|C_t\cup D_w\|>2w-\sqrt{2w}$. In order to prove our main theorem, we generalize this claim.

\begin{prop}\label{lem3}
Suppose the strategy $\sz(w)$ is played resulting in a poset with rainbow chains $C_1,\ldots,C_w$. let $D$ denote a set of $w$ distinct colors, and let $k$ be a positive integer. If $k<\sqrt{2w}$, then there exists a $t$ such that $w-k+1\leq t \leq w$ and $|\sigma(C_t)\cup D|\geq2w-\frac{w}{k}-\frac{k-1}{2}$. Otherwise, there exists an integer $t$ such that $|\sigma(C_t)\cup D|>2w-\sqrt{2w}$.
\end{prop}

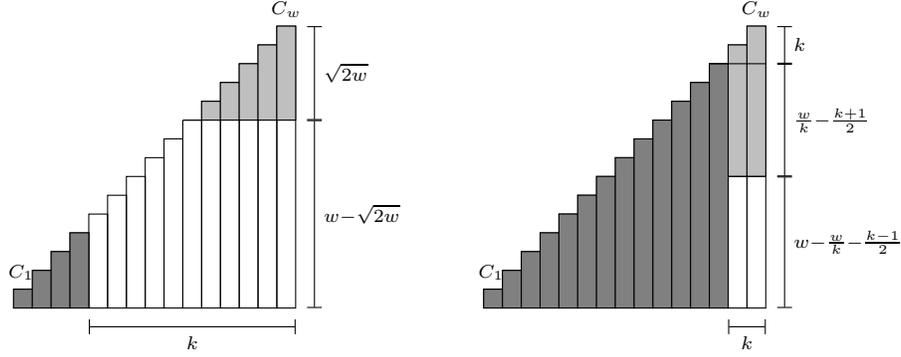
\begin{figure}
\centering
\begin{tikzpicture}[scale=0.25]

\newcommand{\shift}{25}

\foreach \i in {1,...,15}{
    \draw (\i,0) rectangle (\i +1,\i);
    }
\node at (1.4,1) [anchor=south] {$\scriptstyle C_1$} ;
\node at (15.5,15) [anchor=south] {$\scriptstyle C_w$} ;
   
\draw (10,10) -- (16,10);
\foreach \i in {10,...,15}{
    \draw[fill=lightgray] (\i,10) rectangle (\i +1,\i);
    }
\draw[|-|] (17,0) -- (17,10);
\draw (17,5) node[right]{$\scriptstyle w-\sqrt{2w}$};
\draw[|-|] (17,10) -- (17,15);
\draw (17,12.5) node[right]{$\scriptstyle \sqrt{2w}$};

\foreach \i in {1,...,15}{
    \draw (\i + \shift,0) rectangle (\i + \shift + 1,\i);
    }
\node at (1.4+\shift,1) [anchor=south] {$\scriptstyle C_1$} ;
\node at (15.5+\shift,15) [anchor=south] {$\scriptstyle C_w$} ;

\foreach \i in {1,...,13}{
    \draw[fill=gray] (\i + \shift,0) rectangle (\i + \shift +1,\i);
    }
\foreach \i in {1,...,4}{
    \draw[fill=gray] (\i ,0) rectangle (\i +1,\i);
    }

\draw[fill=lightgray] (14+\shift,7) rectangle (15+\shift,14);
\draw[fill=lightgray] (15+\shift,7) rectangle (16+\shift,15);
\draw (13+\shift,13) -- (16+\shift,13);
\draw (14+\shift,7) -- (16+\shift,7);
\draw[|-|] (17+\shift,0) -- (17+\shift,7);
\draw (17+\shift,3.5) node[right]{$\scriptstyle w-\frac{w}{k}-\frac{k-1}{2}$};
\draw[|-|] (17+\shift,13) -- (17+\shift,15);
\draw (17+\shift,14) node[right]{$\scriptstyle k$};
\draw[|-|] (17+\shift,7) -- (17+\shift,13);
\draw (17+\shift,10) node[right]{$\scriptstyle \frac{w}{k} -\frac{k+1}{2}$};
\draw[|-|] (5,-1) -- (16,-1);
\node at (10.5,-1) [anchor=north]{$\scriptstyle k$};
\draw[|-|] (14+\shift,-1) -- (16+\shift,-1);
\node at (15+\shift,-1) [anchor=north]{$\scriptstyle k$};
\end{tikzpicture}
\caption{If $k\geq\sqrt{2w}$, Beth can guarantee the same number of colors as in $\sz'(w)$. If $k<\sqrt{2w}$, Beth can guarantee slightly less colors than $\sz'(w)$ }
\label{fig:n2}
\end{figure}

\begin{proof}
Figure \ref{fig:n2} (LEFT) shows that when $k\geq\sqrt{2w}$, we have enough rainbow chains to guantaree the same number of colors as in $\sz'(w)$. Suppose that $k<\sqrt{2w}$, and let $C'=\bigcup_{i=w-k+1}^w C_i$. Each color from $D$ may only be used once in $C'$ and $|C'|=wk-\frac{1}{2}k(k-1)$. If we let $C''$ denote the set of points not colored with colors from $D$, then $|C''|\geq wk-\frac{1}{2}k(k-1)-w$. On average each chain has $w-\frac{1}{2}(k-1)-\frac{w}{k}$ colors distinct from those in $D$. Thus there must exist an integer $t$ such that $w-k+1\leq t \leq w$ and $|\sigma(C_t)\cup D|\geq2w-\frac{w}{k}-\frac{k-1}{2}$. A visual arguement is shown in Figure \ref{fig:n2} (RIGHT).

\end{proof}

In Stage 3, $\sz'(w)$ ends with Beth recurvisely playing $\sz'(w-1)$ so that every new point is completely incomparable to $C_w\cup D_w$ and completely comparable to $(S_1\cup S_2)\setminus (C_w\cup D_w)$. Any poset resulting from $\sz'(w)$ has width $w$ which is easy to see from the structure of the poset which consists of rainbow chains and the corresponding antichains. We include the following proposition for emphasis.

\begin{prop}\label{prop2}
The strategy $\sz'(w)$ always constructs a poset of width $w$.
\end{prop}
\section{Proof of the First Theorem}

We imitate the improved strategy $\sz'(w)$ from \cite{bos-fel-12} using $L_\alpha(k,w)$ and $L\beta(k,w)$ to obtain a new strategy $S(w)$ for Beth which will force Anna to use $(2-o(1))\binom{w+1}{2}$ colors on a $2$-dimensional poset $(X,P)$ of width $w$. Note that while it may not be necessary to construct the poset via its realizer in order to prove our first result, it is for our main result and hence in using the same substrategies for both results allows us to contrast and better demonstrate the connections between the two problems. 

\subsection{The Strategy for Beth}

We define the strategy $S(w)$ for Beth recursively on the positive integer $w$. The strategy $S(w)$ is completed in three stages. During the first two stages, Beth constructs linear orders $A_k$ and $B_k$ for each $k\in\{1,\ldots,w\}$ such that each pair realizes the same poset $(X,P)$, and presents $(X,P)$ to Anna. In Stage $3$, Beth throws away all but one pair $\{A_t,B_t\}$ of linear orders and plays $S(w-1)$ analogously to that of $\sz'(w)$. The choice of $t$ is dependent on the coloring by Anna. Let $w>N$ for some sufficiently large $N$.

\subsubsection*{Stage 1} 
For each positive integer $k\leq w$, Beth constructs two linear orders $A_k$ and $B_k$ by playing $L_\alpha(k,w)$ and $L_\beta(k,w)$ respectively. By Lemma \ref{lem2}, $S_1$ contains a sequence of rainbow chains $C_1,\ldots,C_w$ such that $C_k<S_1\setminus C_k$ in $A_k$ for $k\in\{1,\ldots,w\}$.

\subsubsection*{Stage 2} 
For every positive integer $k\leq w$, Beth updates $A_k$ and $B_k$ by playing the dual strategies $L_\beta^*(w,w)$ and $L_\alpha^*(w,w)$ completely under $S_1$ in $A_k$ and $B_k$ respectively so that $S_2<S_1$ in both $A_k$ and $B_k$. By Lemma \ref{lem2}, $S_2$ contains a sequence of rainbow chains $D_1,\ldots,D_w$ such that $S_2\setminus D_w<D_w<S_1$ in $B_k$ for $k\in\{1,\ldots,w\}$.

\subsubsection*{Stage 3}
By Proposition \ref{lem3}, there exists a $t$ such that $\|C_t\cup D_w\|>2w-\sqrt{2w}$. Beth plays $S(w-1)$ for the remainder of the game in such a way that $S_2\setminus D_w<S_3<S_1\setminus C_t$ but $S_3$ and $C_t\cup D_w$ are completely incomparable in $P$. 

\subsection{The Result}

Notice that the only difference between $S(w)$ and $\sz'(w)$ is that we kept track of a realizer for each choice of $t$. Thus $S(w)$ forces at least 
\[
\sum_{i=1}^w \left(2w-\sqrt{2w}\right)=\left(2-o(1)\right)\binom{w+1}{2}
\]
colors on a poset $(X,P)$ and by Proposition \ref{prop2}, $(X,P)$ has width at most $w$. 

Finally, we claim that $(X,P)$ is $2$-dimensional. Notice that $A_t\cap B_t=P|_{S_1\cup S_2}$. By the induction hypothesis, the poset $(S_3,P|_{S_3})$ is $2$-dimensional. Let $A$ and $B$ be linear extensions of $P|_{S_3}$ such that $A\cap B = P|_{S_3}$. We define a realizer $\{L_1,L_2\}$ in such a way that $A_t\cup A\subset L_1$, $B_t\cup B\subset L_2$, and the following conditions hold.
\begin{enumerate}
	\item $S_2 < C_t < S_3 < S_1\setminus C_t$ in $L_1$.
	\item $S_2\setminus D_w < S_3 < D_w < S_1$ in $L_2$.
\end{enumerate}

This concludes the proof.

\section{Proof of the Main Theorem}

In this variant of the game, Beth is restricted on the number of linear extensions and must present them to Anna each round. In other words, Beth cannot simply throw away the extra linear extensions while keeping only the ones needed. Hence, we must be more selective when constructing the linear extensions. Similarly to the previous proof, we use the strategies from Section 2 to obtain a new strategy $S(d,w)$ for Beth which will force Anna to use $(2-\frac{1}{d-1}-o(1))\binom{w+1}{2}$ colors on an $d$-dimensional poset $(X,P)$ of width $w$ presented with via a realizer of size $d$.

\subsection{The Strategy for Beth}
We fix the positive integer $d$ and define the strategy $S(d,w)$ for Beth recursively on the positive integer $w$. Beth constructs a poset $(X,P)$ by presenting a realizer $R$ of size $d$. Let $d$ and $w$ be positive integers. Beth constructs $R$ by constructing $d$ linear extensions $L_{w-d+2},\ldots,L_{w},L_{w+1}$. In order to handle the case when $w<d-1$, we extend the algorithm $L_\alpha(k,w)$ to be defined for $k<1$ as follows: If $k<1$, then $L_\alpha(k,w)=L_\alpha(w,w)$. The strategy $S(d,w)$ is completed in three stages.

\subsubsection*{Stage 1} 
For each integer $i$ such that $w-d+2\leq i \leq w$, Beth constructs the linear extension $L_i$ by playing $L_\alpha(i,w)$. Beth simultaneously constructs $L_{w+1}$ by playing $L_\beta(w,w)$. By Lemma \ref{lem2}, $S_1$ contains a sequence of rainbow chains $C_1,\ldots,C_w$ such that $C_i<S_1\setminus C_i$ in $L_i$ for $i\in\{w-d+2,\ldots, w\}$.

\subsubsection*{Stage 2} 
For each integer $i$ such that $w-d+2\leq i \leq w$, Beth updates $L_i$ by playing the dual strategy $L_\beta^*(w,w)$ completely under $S_1$ in $L_i$. Beth simultaneously updates $L_{w+1}$ by playing the dual strategy $L_\alpha^*(w,w)$ completely under $S_1$ in $L_{w+1}$. By Lemma \ref{lem2}, $S_2$ contains a sequence of rainbow chains $D_1,\ldots,D_w$ such that $S_2\setminus D_w<D_w<S_1$ in $L_{w+1}$.

\subsubsection*{Stage 3}

By Proposition \ref{lem3}, for $w$ sufficiently large, there exists a $t$ such that $w-d+2\leq t \leq w$ and $\|C_t\cup D_w\|\geq2w-\frac{w}{d-1}-\frac{d-2}{2}$. Then for each integer $i$ such that $w-d+2\leq i \leq w+1$, Beth plays $L_\alpha(i-1,w-1)$ for each integer $i$ such that $w-d+2\leq i<w$, and $L_\beta(w-1,w-1)$ for $i=w+1$ so that the following inequalities hold:
\begin{enumerate}
    \item $S_2 < C_t < S_3 < S_1\setminus C_t$ in $L_t$

    \item $S_2\setminus D_w < S_3 < D_w < S_1$ in $L_{w+1}$

    \item $S_2 < S_3 < S_1$ in $L_i$ for $i\notin \{t,w+1\}$.
\end{enumerate}

\subsection{The Result}

Since $D_w<C_t<S_3$ in $L_t$ and $S_3<D_w<C_t$ in $L_{w+1}$, $S_3$ and $C_t\cup D_w$ are completely incomparable. Since $S_2\setminus D_w<S_3<S_1\setminus C_t$ in every linear extension, it also holds true in $P$. Hence, the resulting poset $(X,P)$ is in the class $\mathcal{P}$ of posets resulting from $\sz'(w)$. In particular, the posets which can be constructed from $S(d,w)$ are exactly the posets from $\mathcal{P}$ where $t$ is restricted to $\{w-d+2,\ldots,w\}$ when $d<1+\sqrt{2w}$. Thus by Proposition \ref{prop2}, $(X,P)$ has width at most $w$. Moreover, $S(d,w)$ forces at least
\[
\sum_{i=1}^w \left(2i-\frac{i}{d-1}-\frac{d-2}{2}\right)=\left(2-\frac{1}{d-1}-o(1)\right)\binom{w+1}{2}
\]
colors on a poset $(X,P)$ of width at most $w$.

\bibliographystyle{acm}
\bibliography{dd}

\end{document}